\documentclass[11pt]{amsart}
\usepackage{amssymb,latexsym,amsmath,amsfonts,color}
\newtheorem{proposition}{Proposition}[section]
\newtheorem{theorem}[proposition]{Theorem}
\newtheorem{corollary}[proposition]{Corollary}

\newtheorem{remark}[proposition]{Remark}

\newcounter{z}
\newenvironment{enum}{\setcounter{z}{0}
\begin{list}{(\arabic{z})}{\usecounter{z}
\setlength{\topsep}{1ex}\setlength{\labelwidth}{0.5cm}
\setlength{\leftmargin}{1.25cm}\setlength{\labelsep}{0.25cm}
\setlength{\parsep}{2pt}}}{\end{list}~\\[-6ex]}

\newcommand{\bsm}{\left(\begin{smallmatrix}}
\newcommand{\esm}{\end{smallmatrix}\right)}
\newcommand{\beas}{\begin{eqnarray*}}
\newcommand{\eeas}{\end{eqnarray*}}
\newcommand{\Isim}{\stackrel{I}{\sim}}
\newcommand{\s}[1]{\bar{#1}}
\newcommand{\rb}[1]{\raisebox{1.5ex}[-1.5ex]{#1}}

\title{New equivalences for pattern avoiding involutions}

\author{W. M. B. Dukes, V\'{i}t Jel\'{i}nek, Toufik Mansour and Astrid Reifegerste}
\thanks{The second author was supported by project 201/05/H014 of the Czech Science Foundation and project MSM0021620838 of the Czech Ministry of Education.}
\address{Science Institute, University of Iceland, Reykjav\'{i}k, Iceland}
\email{dukes@raunvis.hi.is}
\address{Department of Applied Mathematics, Charles University Prague, Czech Republic}
\email{jelinek@kam.mff.cuni.cz}
\address{Department of Mathematics, University of Haifa, 31905 Haifa, Israel}
\email{toufik@math.haifa.ac.il}
\address{Faculty of Mathematics, University of Magdeburg, Germany}
\email{astrid.reifegerste@ovgu.de}



\begin{document}
\maketitle

\begin{abstract}
We complete the Wilf classification of signed patterns of length 5
for both signed permutations and signed involutions. New general
equivalences of patterns are given which prove Jaggard's
conjectures concerning involutions in the symmetric group avoiding
certain patterns of length 5 and 6. In this way, we also complete
the Wilf classification of $S_5$, $S_6$, and $S_7$ for involutions.
\end{abstract}


\section{Introduction}
Pattern avoidance has proved to be a useful concept in a variety of
seemingly unrelated problems, including Kazhdan-Lusztig 
polynomials~\cite{Bec1997}, singularities of Schubert varieties
\cite{Bil1998,BilJocSta1993,BilKai1988,BilLak,BilWar2000,LakSan1990}, 
Chebyshev polynomials \cite{ManVai2002b}, rook polynomials
for a rectangular board \cite{ManVai2000a} and various 
sorting algorithms, sorting stacks and sortable permutations 
\cite{Bon2002a,Bon2002c,Bou1998,Tar1972,Wes1990b,Wes1993}.

In this paper, we deal with pattern avoidance in the symmetric group $S_n$ and
the hyperoctahedral group $B_n$.
The group $B_n$, which is isomorphic to the automorphism group of the
$n$-dimensional hypercube, can be represented as the group of all bijections $\omega$ of the set 
$X=\{-n,\ldots , -1,1,\ldots,n\}$ onto itself such that 
$\omega(-i)\,=\,-\omega(i)$ for all $i \in X$, with composition 
as the group operation.
However, for our purposes it is more convenient to represent the elements of $S_n$ as permutation matrices, and the
elements of $B_n$ as signed permutation matrices, where a signed
permutation matrix is a ${0,1,-1}$-matrix with exactly one nonzero
entry in every row and every column. We may also write the
elements of $B_n$ as words $\pi=\pi_1\pi_2\ldots\pi_n$ in which
each of the letters $1,2,\ldots,n$ appears, possibly barred to
signify negative letters; a matrix $p$ corresponds to the word
$\pi$ such that $p_{ij}=1$ if $\pi_i=j$, $p_{ij}=-1$ if
$\pi_i=-j$, and $p_{ij}=0$ otherwise. In our paper, we will make
no explicit distinction between these two representations of a
signed permutation. Let $I_n$ and $SI_n$ be the set of involutions
in $S_n$ and $B_n$, respectively. Note that involutions correspond
precisely to symmetric matrices.

A signed permutation $\pi\in B_n$ is said to {\it contain the
pattern} $\tau\in B_k$ if there exists a sequence $1\le
i_1<i_2<\ldots<i_k\le n$ such that $|\pi_{i_a}|<|\pi_{i_b}|$ if
and only if $|\tau_a|<|\tau_b|$ and $\pi_{i_a}>0$ if and only if
$\tau_a>0$ for all $1\le a,b\le k$. Otherwise, $\pi$ is called a
{\it $\tau$-avoiding} permutation. Note that $\pi$ contains $\tau$
if and only if the matrix representing $\pi$ contains the matrix
representing $\tau$ as a submatrix. By $M(\tau)$ we denote the set
of all elements of $M$ which avoid the pattern $\tau$.

Two signed patterns $\sigma$ and $\tau$ are called {\it Wilf
equivalent}, in symbols $\sigma\sim\tau$, if they are avoided by
the same number of signed $n$-permutations, i.e., if
$|B_n(\sigma)|=|B_n(\tau)|$ for each $n\ge 1$. Similarly, $\sigma$
and $\tau$ are called \emph{I-Wilf equivalent}, denoted by $\sigma
\Isim\tau$, if $|SI_n(\sigma)|=|SI_n(\tau)|$ for each $n$. Note
that two unsigned permutations $\sigma,\tau\in S_k$ are
Wilf-equivalent if and only if they satisfy the identity
$|S_n(\sigma)|=|S_n(\tau)|$ for each $n$, and they are I-Wilf
equivalent if and only if they satisfy $|I_n(\sigma)|=|I_n(\tau)|$
for each $n$. The classification given by the Wilf equivalence is
slightly coarser than that which is based on the symmetries of
permutations, that is, the mappings generated by the reversal,
transpose, and barring operation. The same is true for the I-Wilf
equivalence, where the available symmetries are generated by the
two diagonal reflections and the barring operation.

The question of whether two patterns are Wilf equivalent or not is
difficult to answer in many cases. By the few generic equivalences
known so far, it has been possible to completely determine the
Wilf classes of $S_n$ up to level $n=7$. The decomposition of
$S_n$ into I-Wilf classes has been completely determined for $n=4$
and almost solved for $n=5$ as well. Jaggard \cite{jaggard}
conjectured the last case of a possible equivalence for patterns
of length 5: $12345$ (or equivalently, $54321$) and $45312$ are
equally restrictive for $I_n$ up to $n=11$.

Continuing the I-Wilf classification of signed patterns that began
in \cite{dukes-mansour-reifegerste}, we will first prove a general
equivalence result which confirms Jaggard's conjecture mentioned
above, as well as another conjecture he made about the equivalence
of certain patterns of length 6. The correspondence behind this
result is based on a bijection between pattern avoiding
transversals of Young diagrams given by Backelin, West and Xin
\cite{backelin-west-xin}. In this way, we complete the
classification of $S_5$ with respect to $\Isim$, which is
fundamental for the analogous classification of $B_5$. The result
even covers all missing I-Wilf equivalences in $S_6$ and $S_7$.

Furthermore, we will show that barring some blocks of a signed
block diagonal pattern preserves the Wilf class of the pattern,
and it also (under some additional assumptions) preserves the
I-Wilf class. These results not only allow us to determine the
Wilf as well as the I-Wilf classes in $B_5$ but they also have
consequences for longer signed patterns.


\section{Jaggard's conjectures}

In 2003, Jaggard \cite{jaggard} proved the equivalences
$12\tau\Isim 21\tau$ and $123\tau\Isim 321\tau$, and completed the
classification of $S_4$ according to pattern avoidance by
involutions in this way. Furthermore, he conjectured that
\vspace*{-1ex}

\begin{enum}
\item $12\ldots k\tau\Isim k(k-1)\ldots1\tau$ for any $k\ge 1$,
\item $12345\Isim 45312$ (or equivalently, $54321\Isim 45312$),
\item $123456\Isim 456123\Isim 564312$ (or equivalently, $654321\Isim 456123$).
\end{enum}
\vspace*{-2ex}

In \cite{backelin-west-xin}, Backelin, West and Xin defined a
transformation to prove $12\ldots k\tau\sim k(k-1)\ldots1\tau$.
(As already mentioned in \cite{dukes-mansour-reifegerste}, their
proof also works for a signed pattern $\tau$.) This map acts not
only on permutation matrices, but more generally, on transversals
of Young diagrams. Bousquet-M{\'e}lou and Steingr{\'i}msson
\cite{bousquet-steingrimsson} showed that this map commutes with
the diagonal reflection of the diagram, which proves the first of
the three conjectures above. From this result, it follows that
$$
\bsm\alpha_k&0&0\\0&\chi&0\\0&0&\alpha_l\esm\Isim\bsm\beta_k&0&0\\0&\chi&0\\0&0&\beta_l\esm
$$
for every signed permutation matrix $\chi$ and any $k,l\ge0$,
where $\alpha_n$ and $\beta_n$ denote the $n\times n$ diagonal and
antidiagonal permutation matrices corresponding to $12\ldots n$
and $n(n-1)\ldots 1$, respectively. In this section, we will show
that
$$
\bsm0&0&0&\alpha_k\\0&0&\chi&0\\0&\chi^t&0&0\\\alpha_k&0&0&0\esm\Isim
\bsm0&0&0&\beta_k\\0&0&\chi&0\\0&\chi^t&0&0\\\beta_k&0&0&0\esm\quad\mbox{ and }\quad
\bsm0&0&0&0&\alpha_k\\0&0&0&\chi&0\\0&0&1&0&0\\0&\chi^t&0&0&0\\\alpha_k&0&0&0&0\esm\Isim
\bsm0&0&0&0&\beta_k\\0&0&0&\chi&0\\0&0&1&0&0\\0&\chi^t&0&0&0\\\beta_k&0&0&0&0\esm ,
$$
where $\chi^t$ denotes the transpose of $\chi$. Note that,
different to the general case, the reverse operation is not a
symmetry for involutions, so these equivalences are really new.

Our proof will also use the Backelin, West and Xin bijection~\cite{backelin-west-xin}. 
Therefore, let us first
recall the extended notion of pattern avoidance they have used. A
{\it Young diagram} (or Young shape) is a top-justified and
left-justified array of cells, i.e., an array whose rows have
non-increasing lengths from top to bottom, and its columns have
non-increasing lengths from left to right. A cell of a Young shape
is called a {\it corner} if the array obtained by removing the
cell is still a Young shape. Occasionally, it will be convenient
to use top-right justified diagrams instead of the top-left
justified diagrams defined above. We will refer to the top-right
justified shapes as {\it NE-shapes} to avoid confusion with the
ordinary Young shapes.

A ({\it signed}) {\it transversal} of a Young diagram $\lambda$ is
an assignment of 0's and 1's (of 0's, 1's and -1's) to the cells
of $\lambda$, such that each row and column contains exactly one
nonzero entry. A \emph{sparse filling} of $\lambda$ is an
arrangement of 0's, 1's and -1's which has at most one nonzero
entry in every row and column.

For a $k\times k$ permutation matrix $\tau$, we say that a filling
$L$ of a shape $\lambda$ {\it contains $\tau$} if there exists a
$k\times k$ subshape within $\lambda$ whose induced filling is
equal to $\tau$. The set of all transversals (or signed
transversals) of a shape $\lambda$ which do not contain $\tau$ is
denoted by $S_\lambda(\tau)$ (or $B_\lambda(\tau)$, respectively).
Two signed permutation matrices $\sigma$ and $\tau$ are called
{\it shape Wilf equivalent} if
$|B_\lambda(\sigma)|=|B_\lambda(\tau)|$ for all Young shapes
$\lambda$. Shape Wilf equivalence clearly implies Wilf
equivalence. We will also say that $\sigma$ and $\tau$ are
\emph{NE-shape Wilf equivalent} if
$|B_\lambda(\sigma)|=|B_\lambda(\tau)|$ for each NE-shape
$\lambda$. Observe that if $\sigma$ and $\tau$ are permutation
matrices, then they are shape Wilf equivalent if and only if
$|S_\lambda(\sigma)|=|S_\lambda(\tau)|$ for each Young diagram
$\lambda$.

By \cite[Proposition 2.2]{backelin-west-xin}, $\alpha_k$ and
$\beta_k$ are shape Wilf equivalent for all $k$. The following
proposition, which is also largely based on
\cite{backelin-west-xin}, will allow us to extend this equivalence
to more general patterns.

\begin{proposition} \label{bwx}
Let $\lambda$ be a Young shape, and let $\chi,\chi_1,\chi_2$ be
signed permutations, such that $\chi_1$ and $\chi_2$ are shape
Wilf equivalent. We set
$$
\theta=\bsm\chi_1&0\\0&\chi\esm\quad\mbox{and}\quad\omega=\bsm\chi_2&0\\0&\chi\esm.
$$
There is a bijection between $\theta$-avoiding and
$\omega$-avoiding sparse fillings of $\lambda$. This bijection
preserves the number of nonzero entries in each row and column; in
particular, $\theta$ and $\omega$ are shape Wilf equivalent.
Furthermore, if $\chi$ is nonempty, the bijection preserves the
values of the filling in the corners of $\lambda$.
\end{proposition}

\begin{proof}
The proof is essentially the same as the proof given in
\cite[Proposition 2.3]{backelin-west-xin}. We briefly sketch the
argument here. By assumption, there is a bijection $\varphi$
between the $\chi_1$-avoiding and $\chi_2$-avoiding signed
transversals of an arbitrary Young shape. Let $L$ be an arbitrary
$\theta$-avoiding sparse filling of $\lambda$. Let us colour a
cell of $\lambda$ if there is no occurrence of $\chi$ to the
south-east of this cell. Also, if $\lambda$ has a row or column
where all the uncoloured cells contain zeros, then we colour each
cell of this row or column. Note that if $\chi$ is nonempty, then
all the corners of $\lambda$ are coloured. The uncoloured cells
induce a $\chi_1$-avoiding signed transversal of a Young
subdiagram of $\lambda$. We apply the bijection $\varphi$ to the
subdiagram of uncoloured cells, and preserve the filling of all
the coloured cells. This transforms the original filling of
$\lambda$ into a $\omega$-avoiding sparse filling. This
transformation is a bijection which has all the claimed
properties.
\end{proof}

Note that Proposition~\ref{bwx} yields some information even when
$\chi$ is the empty matrix. In such situation, the proposition
shows that a bijection between pattern avoiding signed
transversals can be extended to a bijection between
pattern-avoiding sparse fillings, by simply ignoring the rows and
columns with no nonzero entries.

We will now show how the results on shape Wilf equivalence may be
applied to obtain new classes of I-Wilf equivalent patterns. Let
us first give the necessary definitions. For an $n\times n$ matrix
$\pi$ let $\pi^+$ denote the subfilling of $\pi$ formed by the
cells of $\pi$ which are strictly above the main diagonal, and let
$\pi_0^+$ denote the subfilling formed by the cells on the main
diagonal and above it. For example, for $\pi=2\bar{4}31$ we have
\begin{center}
\unitlength=0.4cm
\begin{picture}(16,3)
\put(2,3){\line(1,0){3}}\put(2,2){\line(1,0){3}}\put(3,1){\line(1,0){2}}\put(4,0){\line(1,0){1}}
\put(2,2){\line(0,1){1}}\put(3,1){\line(0,1){2}}\put(4,0){\line(0,1){3}}\put(5,0){\line(0,1){3}}
\put(0,1.5){\makebox(0,0)[cc]{$\pi^+=$}}
\put(2.5,2.5){\makebox(0,0)[cc]{\footnotesize $1$}}
\put(4.5,1.5){\makebox(0,0)[cc]{\footnotesize $-1$}}
\put(13,3){\line(1,0){4}}\put(13,2){\line(1,0){4}}\put(14,1){\line(1,0){3}}\put(15,0){\line(1,0){2}}\put(16,-1){\line(1,0){1}}
\put(13,2){\line(0,1){1}}\put(14,1){\line(0,1){2}}\put(15,0){\line(0,1){3}}\put(16,-1){\line(0,1){4}}\put(17,-1){\line(0,1){4}}
\put(9.5,1.5){\makebox(0,0)[cc]{$\mbox{ and }\quad\pi_0^+=$}}
\put(14.5,2.5){\makebox(0,0)[cc]{\footnotesize $1$}}
\put(16.5,1.5){\makebox(0,0)[cc]{\footnotesize $-1$}}
\put(15.5,0.5){\makebox(0,0)[cc]{\footnotesize $1$}}
\put(17.5,1.3){\makebox(0,0)[cc]{$\mbox{.}$}}
\end{picture}
\end{center}
\vspace*{2ex} The coordinates of the entries in $\pi$ are used for
the cells of $\pi^+$ as well. Thus, for instance, the cell $(1,2)$
is the top-left corner of $\pi^+$. Analogously, we define $\pi^-$
to be the filled shape corresponding to the entries strictly below
the main diagonal of $\pi$. Clearly, a symmetric matrix $\pi$ is
completely determined by $\pi_0^+$. Observe that a symmetric
${0,1,-1}$-matrix $\pi$ is a signed involution if and only if, for
every $i=1,\ldots,n$, the filling $\pi_0^+$ has exactly one
nonzero entry in the union of all cells of the $i$-th row and
$i$-th column.

Note that $i$ is a fixed point of a signed involution $\pi$, that
is $|\pi_i|=i$, if and only if the $i$-th row and the $i$-th
column of $\pi^+$ have all entries equal to zero. In general, a
signed involution $\pi$ need not be completely determined by the
filling $\pi^+$; however, if we have two signed involutions
$\pi,\rho$ with $\pi^+=\rho^+$, then $\pi$ and $\rho$ only differ
by the signs of their fixed points. If $\pi$ is a signed
involution, then, for each $i=1,\dotsc,n$, the filling $\pi^+$ has
at most one nonzero entry in the union of the $i$-th row and
$i$-th column; conversely, any filling $\pi^+$ of appropriate
shape with these properties can be extended into a signed
involution $\pi$, which is determined uniquely up to the sign of
its fixed points.

For a signed permutation $\sigma$, let $\sigma'$ denote the
involution $\bsm 0 &\sigma\\ \sigma^t & 0\esm$, where $\sigma^t$
is the transpose of $\sigma$. We are now ready to state our first
result on I-Wilf equivalence.

\begin{theorem} \label{equivalence1}
If $\sigma$ and $\tau$ are two NE-shape Wilf equivalent signed
permutation matrices, then $\sigma'\Isim\tau'$. Moreover, the
bijection between $SI_n(\sigma')$ and $SI_n(\tau')$ preserves
fixed points.
\end{theorem}

\begin{proof}
Let $\pi\in SI_n$ be an involution. We claim that $\pi$ avoids
$\sigma'$ if and only if $\pi^+$ avoids $\sigma$. To see this,
notice that any occurrence of $\sigma'$ in $\pi$ can be restricted
either to an occurrence of $\sigma$ in $\pi^+$ or an occurrence of
$\sigma^t$ in $\pi^-$; however, since $\pi^+$ is the transpose of
$\pi^-$, we know that $\pi^-$ contains $\sigma^t$ if and only if
$\pi^+$ contains $\sigma$. The converse is even easier to see.

Let us choose $\pi\in SI_n(\sigma')$. Since $\pi^+$ is a sparse
$\sigma$-avoiding filling, we may apply the bijection from
Proposition \ref{bwx} (adapted for NE-shapes) to $\pi^+$, to
obtain a $\tau$-avoiding sparse filling of the same shape, which
has a nonzero entry in a row $i$ (or column $i$) whenever $\pi^+$
has a nonzero entry in the same row (or column, respectively).
Hence this filling also corresponds to an involution, more
exactly, to $\rho^+$ for an involution $\rho\in SI_n$, and
furthermore, the fixed points of $\rho$ are in the same position
as the fixed points of $\pi$, because the position of the fixed
points is determined by the zero rows and columns, which are
preserved by the bijection from Proposition~\ref{bwx}. By defining
the signs of the fixed points of $\rho$ to be the same as the
signs of the fixed points of $\pi$, the involution $\rho$ is
determined uniquely. Clearly, since $\rho^+$ avoids $\tau$, we
know that $\rho$ avoids $\tau'$. Each step of this construction
can be inverted which proves the bijectivity. Furthermore, the
bijection preserves fixed points by construction.
\end{proof}

By a similar reasoning, we obtain an analogous result for patterns
of odd size. For a signed permutation $\sigma$, let $\sigma''$
denote the involution matrix
$$
\bsm 0 & 0 & \sigma \\ 0 & 1 & 0 \\  \sigma^t &0&0\esm,
$$
and let $\sigma^*$ denote the signed permutation $\bsm 0 &\sigma
\\ 1 & 0\esm$.

\begin{theorem} \label{equivalence2}
If $\sigma$ and $\tau$ are NE-shape Wilf equivalent, then
$\sigma''\Isim\tau''$. Moreover, the bijection between
$SI_n(\sigma'')$ and $SI_n(\tau'')$ preserves fixed points.
\end{theorem}
\begin{proof}
By an argument analogous to the proof of
Theorem~\ref{equivalence1}, we may observe that an involution
$\pi$ avoids $\sigma''$ if and only if $\pi_0^+$ avoids the
pattern $\sigma^*$. By Proposition~\ref{bwx} (adapted for
NE-shapes), the two patterns $\sigma^*$ and $\tau^*$ are NE-shape
Wilf equivalent and furthermore, the bijection realizing this
equivalence preserves the corners of the shape. Note that in our
situation, the corners correspond exactly to the diagonal cells of
the original signed permutation matrix.

Now we consider $\pi_0^+$ for an involution $\pi\in
SI_n(\sigma'')$. By Proposition \ref{bwx}, $\pi_0^+$ is in
bijection with a $\tau^*$-avoiding filling $\rho_0^+$. Since the
bijection preserves the number of nonzero entries in each row and
each column of $\pi^+_0$, and it also preserves the entries on the
intersection of $i$-th row and $i$-th column (these are precisely
the corners), we know that the bijection preserves, for each $i$,
the number of nonzero entries in the union of the $i$-th row and
$i$-th column. In particular, $\rho_0^+$ has exactly one nonzero
entry in the union of $i$-th row and $i$-th column, which
guarantees that $\rho_0^+$ can be (uniquely) extended into an
involution $\rho$.

Because the bijection preserves the entries in the diagonal cells
$(i,i)$, $i=1,\ldots,n$, the permutations $\pi$ and $\rho$ have
the same fixed points. This provides the required bijection.
\end{proof}

Let us apply these two theorems to some special cases of shape
Wilf equivalent patterns. For an integer $k\ge0$ and a signed
permutation $\chi$, let us define
$$
\theta=\bsm0&\alpha_k\\\chi&0\esm\quad\mbox{and}\quad\omega=\bsm0&\beta_k\\\chi&0\esm.
$$
As we know, the two patterns $\theta$ and $\omega$ are NE-shape
Wilf equivalent. From our results, we then obtain the following
classes of I-Wilf equivalent patterns.

\begin{corollary}We have
$$
\bsm0&0&0&\alpha_k\\0&0&\chi&0\\0&\chi^t&0&0\\\alpha_k&0&0&0\esm\Isim
\bsm0&0&0&\beta_k\\0&0&\chi&0\\0&\chi^t&0&0\\\beta_k&0&0&0\esm\quad
\mbox{and}\quad
\bsm0&0&0&0&\alpha_k\\0&0&0&\chi&0\\0&0&1&0&0\\0&\chi^t&0&0&0\\\alpha_k&0&0&0&0\esm\Isim
\bsm0&0&0&0&\beta_k\\0&0&0&\chi&0\\0&0&1&0&0\\0&\chi^t&0&0&0\\\beta_k&0&0&0&0\esm.
$$
\end{corollary}

The special cases $\chi=\emptyset$ and $\chi=(1)$ show both of
Jaggard's conjectures to be correct.

\begin{corollary}
We have $54321\Isim 45312$ and $654321\Isim 456123\Isim 564312$.
\end{corollary}


\section{Barring of blocks}

In \cite{dukes-mansour-reifegerste} it was shown that the barring
of $\tau$ in $12\ldots k\tau$ and $k(k-1)\ldots1\tau$ preserves
both the Wilf class and the I-Wilf class. Furthermore it was
proved that
$$
\bsm\alpha_k&0&0\\0&\chi&0\\0&0&\alpha_k\esm\Isim\bsm\alpha_k&0&0\\0&-\chi&0\\0&0&\alpha_k\esm
$$
for every signed permutation matrix $\chi$ and $k\ge0$. Basically,
the assertion follows from $123\Isim1\bar{2}3$. By a similar
reasoning, we can show the I-Wilf equivalence of the reversed
patterns because $321\Isim3\bar{2}1$ as well. Now we turn our
attention to the general block pattern
$$
\bsm\chi_1&0&0\\0&\chi_2&0\\0&0&\chi_3\esm
$$
where the $\chi_i$ are signed permutation matrices. First we prove
the following crucial statement.

\begin{theorem} \label{barring}
Let $\chi_1$ and $\chi_2$ be signed permutation matrices and set
$$
\theta=\bsm\chi_1&0\\0&\chi_2\esm\quad\mbox{and}\quad\omega=\bsm\chi_1&0\\0&-\chi_2\esm.
$$
For any Young shape $\lambda$, there is a bijection between
$\theta$-avoiding and $\omega$-avoiding sparse fillings
of~$\lambda$. The bijection preserves the position of all nonzero
entries, i.e., it transforms the filling only by changing the
signs of some of the entries. In particular, the patterns $\theta$
and $\omega$ are shape Wilf equivalent. Moreover, if $\lambda$ is
self-conjugate and at least one of the matrices $\chi_1$ and
$\chi_2$ is symmetric, then the bijection maps symmetric fillings
to symmetric fillings.
\end{theorem}

\begin{proof}
Given a $\theta$-avoiding sparse filling of $\lambda$, we
construct the corresponding $\omega$-avoiding filling as follows:
Colour each cell of $\lambda$ for which there is an occurrence of
$\chi_1$ to the north-west of the cell. Note that the cells left
uncoloured then form a Young subdiagram of $\lambda$. By
assumption, the coloured part does not contain $\chi_2$. Switching
the signs of all entries of this part consequently yields a signed
transversal of $\lambda$ which avoids $\omega$. Note that even
after the transformation has been performed, it is still true that
the coloured cells are precisely those cells that have an
occurrence of $\chi_1$ to their north-west. The transformation may
have created new copies of $\chi_1$ in the diagram, but it may be
easily seen that these copies do not alter the colouring of the
cells. This shows that the transformation is indeed a bijection.

Let $\lambda$ now be self-conjugate with a symmetric
$\theta$-avoiding filling. Obviously, if $\chi_1$ is symmetric,
then a cell is coloured if and only if its reflection (along the
main diagonal) is coloured. Hence the signs of both entries must
have been changed, so the resulting filling is symmetric again. If
$\chi_2$ is symmetric but $\chi_1$ is not, then we slightly modify
the definition of the bijection. Colour a cell if there is an
occurrence of $\chi_2$ to the south-east. The restriction to these
cells is a symmetric filling of a self-conjugate subshape which
avoids $\chi_1$. Now change the signs of all nonzeros in
uncoloured cells. The resulting filling avoids $\omega$ and is
still symmetric. It is again easy to see that this provides the
required symmetry-preserving bijection.
\end{proof}

An immediate consequence of the previous theorem is the following:

\begin{corollary} \label{equivalence3}
For any signed permutation matrices $\chi_1,\chi_2,\chi_3$, we
have
$$
\bsm\chi_1&0&0\\0&\chi_2&0\\0&0&\chi_3\esm\sim\bsm\chi_1&0&0\\0&-\chi_2&0\\0&0&\chi_3\esm.
$$
\end{corollary}

Because of the symmetry property of the bijection we can prove an
analogous result for pattern avoiding involutions.

\begin{corollary} \label{equivalence4}
Let $\chi_1,\chi_2,\chi_3$ be signed permutation matrices, at least two of which are symmetric. Then we have
$$
\bsm\chi_1&0&0\\0&\chi_2&0\\0&0&\chi_3\esm\Isim\bsm\chi_1&0&0\\0&-\chi_2&0\\0&0&\chi_3\esm.
$$
\end{corollary}

\begin{proof}
By Theorem \ref{barring}, the signed pattern ${\rm
diag}(\chi_1,\chi_2,\chi_3)$ is I-Wilf equivalent with the signed
pattern ${\rm diag}(\chi_1,\chi_2,-\chi_3)$ (note that at least
one of the two matrices ${\rm diag}(\chi_1,\chi_2)$ and $\chi_3$
is symmetric). By the same argument, the pattern ${\rm
diag}(\chi_1,\chi_2,\chi_3)$ is I-Wilf equivalent with ${\rm
diag}(\chi_1,-\chi_2,-\chi_3)$. Combining these facts with the
observation that changing the signs of all the three blocks
clearly preserves the I-Wilf class, we may even conclude that any
matrix obtained by changing the signs of any of the three blocks
is I-Wilf equivalent with the original matrix.
\end{proof}

Combining Theorem~\ref{barring} with Theorems~\ref{equivalence1}
and~\ref{equivalence2}, we obtain more classes of I-Wilf
equivalent patterns. The following corollary gives an example.

\begin{corollary} \label{equivalence5}
Let $\chi_1$ and $\chi_2$ be signed permutation matrices. Then we have
$$
\bsm0&0&0&0&\chi_1\\0&0&0&\chi_2&0\\0&0&\varepsilon&0&0\\0&\chi_2^t&0&0&0\\\chi_1^t&0&0&0&0\esm\Isim\bsm0&0&0&0&\chi_1\\0&0&0&-\chi_2&0\\0&0&\varepsilon&0&0\\0&-\chi_2^t&0&0&0\\\chi_1^t&0&0&0&0\esm.
$$
where $\varepsilon$ is empty or $\varepsilon=(1)$.
\end{corollary}


\section{Classification}

The proof of Jaggard's conjecture provides the complete
classification of the I-Wilf equivalences among the patterns from
$S_5$. It turns out that there are 36 different classes (in
comparison with 45 symmetry classes). By the results of
\cite{dukes-mansour-reifegerste}, it has been known that $B_5$ has
at most 405 I-Wilf equivalence classes. Applying the new
equivalences, we obtain 402 classes which are definitively
different. (By the symmetries of an involutive permutation, the
patterns are divided into 566 classes.) Table 1 shows
representatives of all classes, each with the number of
involutions in $SI_9,\ldots,SI_{12}$ avoiding the patterns of this
class. The enumeration is done for $n=9$ in any case; higher
levels are only computed up to the final distinction. Classes
containing patterns of $S_5$ are in bold; hence the classification
of $S_5$ according to the I-Wilf equivalence can be read from the
table as well.

The classification of the patterns of $B_5$ by Wilf equivalence
becomes complete by Corollary~\ref{equivalence3}. The relations
given in \cite{dukes-mansour-reifegerste} did not cover seven
pairs of patterns whose Wilf equivalence was indicated by
numerical results. All these cases are proved now by the
corollary. Consequently, $B_5$ falls into 130 Wilf classes (in
comparison with 284 symmetry classes). See \cite[Table
7]{dukes-mansour-reifegerste} for the complete list.

The bijections of Theorem \ref{equivalence1} and Theorem
\ref{equivalence2} also provide the complete classification of
$S_6$ and $S_7$ with respect to the I-Wilf equivalence. Table 2
lists all classes of $S_6$ obtained by all equivalences, already
known (see \cite{dukes-mansour-reifegerste} and the references therein) or proven here. As the enumeration of involutions in
$I_{12}$ avoiding the patterns shows, they are different. In a
similar way, we obtain  $1291$ Wilf classes for $S_7$ whose table
is available from \cite{www}.

It is very possible that the results given here and in
\cite{dukes-mansour-reifegerste} suffice to solve the I-Wilf
classification of signed patterns up to length 7. However, the
numerical proof that two classes are really different for a
rapidly increasing number of classes is the challenge we (and
computers) have to master.

\begin{remark}
{\rm After publishing this paper in arXiv, Aaron Jaggard mentioned
that he and Joseph Marincel had shown that the patterns
$(k-1)k(k-2)\ldots312$ and $k(k-1)\ldots21$ are I-Wilf equivalent
for any $k\ge5$ by using generating tree techniques
\cite{jaggard-marincel}.}
\end{remark}

\begin{table}[ht]
{\tiny
\begin{tabular}{|cl|cl|cl|cl|}\hline
$3\s{5}1\s{4}\s{2}$&$160482$&{\bf35142}&$160519$&$14\s{5}2\s{3}$&$160623$&$3\s{5}14\s{2}$&$160627$\\\hline
$4\s{5}31\s{2}$&$160647$&$351\s{4}2$&$160662$&{\bf14325}&$160668$&{\bf12435}&$160670$\\\hline
$52\s{4}\s{3}1$&$160682$&&$160684$ $856400$&&$160684$ $856400$&$523\s{4}1$&$160684$ $856396$\\[-1ex]
&&\rb{\bf 12345}&$4724160$&\rb{\bf 52431}&$4724162$&&\\\hline
$52\s{3}41$&$160686$&{\bf52341}&$160702$&$15\s{3}42$&$160817$&{\bf14523}&$160819$\\\hline
$153\s{4}2$&$160831$&{\bf15342}&$160834$&$125\s{4}3$&$160843$&$15\s{4}\s{3}2$&$160845$\\\hline
$15\s{3}\s{4}2$&$160861$&$14\s{3}25$&$160944$&$124\s{3}5$&$164848$&$13\s{4}25$&$165194$\\\hline
$1325\s{4}$&$165198$&$13\s{5}\s{4}\s{2}$&$165227$&$1235\s{4}$&$165230$&{\bf13542}&$165269$\\\hline
$524\s{3}1$&$165304$&{\bf13425}&$165310$&$124\s{5}\s{3}$&$165365$&$143\s{5}2$&$165389$\\\hline
{\bf14352}&$165416$&$154\s{3}2$&$165484$&$124\s{5}3$&$165525$&$2\s{5}\s{4}\s{3}\s{1}$&$165557$\\\hline
{\bf25431}&$165560$&$13\s{5}2\s{4}$&$165585$&$2\s{5}1\s{4}\s{3}$&$165588$&{\bf45231}&$165596$\\\hline
{\bf12453}&$165598$&$15\s{4}3\s{2}$&$165600$&$2\s{1}5\s{4}3$&$165604$&{\bf25143}&$165627$\\\hline
$4\s{5}2\s{3}1$&$165734$&$53\s{4}21$&$165777$&{\bf13524}&$165788$&{\bf53421}&$165990$\\\hline
$143\s{2}5$&$166106$&$134\s{2}5$&$166279$&$1254\s{3}$&$166337$&$13\s{5}4\s{2}$&$166363$\\\hline
$13\s{4}5\s{2}$&$166398$&$13\s{4}\s{5}\s{2}$&$166404$ $896272$&$135\s{4}2$&$166404$ $896308$&$134\s{5}2$&$166418$\\\hline
{\bf13452}&$166429$&$14\s{5}3\s{2}$&$166451$&$2\s{5}14\s{3}$&$166467$&{\bf14532}&$166479$\\\hline
$3\s{5}2\s{4}\s{1}$&$166488$&$1245\s{3}$&$166498$&$251\s{4}3$&$166505$&$14\s{3}\s{5}2$&$166527$ $897293$\\\hline
{\bf35241}&$166527$ $897923$&$14\s{3}52$&$166538$&$143\s{5}\s{2}$&$166544$&$1543\s{2}$&$166550$\\\hline
$2\s{5}3\s{4}\s{1}$&$166567$&$25\s{3}41$&$166569$&$13\s{5}\s{4}2$&$166572$&{\bf32541}&$166575$\\\hline
$2\s{5}\s{3}4\s{1}$&$166581$&$253\s{4}1$&$166583$&$24\s{5}1\s{3}$&$166586$&$2\s{5}\s{3}\s{4}\s{1}$&$166587$\\\hline
$134\s{5}\s{2}$&$166591$ $898088$&{\bf25341}&$166591$ $898195$&$14\s{5}\s{3}\s{2}$&$166607$&$13\s{4}52$&$166615$\\\hline
$13\s{4}\s{5}2$&$166619$&$14\s{5}\s{3}2$&$166627$&{\bf24513}&$166628$ $898700$&$543\s{2}1$&$166628$ $898668$\\\hline
$145\s{3}2$&$166655$&$3\s{5}24\s{1}$&$166658$&$352\s{4}1$&$166662$&$135\s{2}4$&$166701$\\\hline
$2\s{5}43\s{1}$&$166720$&$25\s{4}\s{3}1$&$166723$&$1354\s{2}$&$166725$ $899209$&$1435\s{2}$&$166725$ $899210$\\\hline
$13\s{5}\s{2}\s{4}$&$166727$&$2\s{5}34\s{1}$&$166737$&$25\s{3}\s{4}1$&$166739$&$2\s{4}5\s{1}3$&$166741$\\\hline
$32\s{5}\s{4}1$&$166742$&$25\s{1}43$&$166754$&$14\s{5}32$&$166755$&$2\s{5}\s{1}\s{4}\s{3}$&$166756$\\\hline
$2\s{4}\s{3}5\s{1}$&$166757$&$243\s{5}1$&$166758$&$2\s{3}\s{5}\s{4}\s{1}$&$166759$ $899733$&{\bf24351}&$166759$ $899753$\\\hline
{\bf23541}&$166760$&$2\s{4}\s{5}\s{1}\s{3}$&$166761$&$2\s{3}5\s{1}4$&$166762$&$23\s{5}1\s{4}$&$166769$\\\hline
$1345\s{2}$&$166773$ $899813$&$2\s{5}4\s{3}\s{1}$&$166773$ $899906$&$25\s{4}31$&$166775$ $899951$&$534\s{2}1$&$166775$ $900042$\\\hline
$23\s{5}\s{4}\s{1}$&$166776$&$2\s{3}541$&$166777$&$2\s{3}\s{5}\s{1}\s{4}$&$166780$&{\bf45321}&$166788$\\\hline
$5432\s{1}$&$166790$&{\bf23514}&$166791$&$4\s{5}32\s{1}$&$166800$&$3\s{5}\s{4}1\s{2}$&$166805$\\\hline
{\bf35412}&$166809$&$2\s{5}\s{4}1\s{3}$&$166816$&$35\s{2}41$&$166818$&{\bf25413}&$166822$\\\hline
$3\s{5}\s{2}\s{4}\s{1}$&$166834$&$25\s{4}13$&$166861$&$13\s{5}24$&$166863$&$1352\s{4}$&$166875$\\\hline
$2\s{5}41\s{3}$&$166876$&$2\s{3}54\s{1}$&$166933$&$23\s{5}\s{4}1$&$166934$ $901415$&$2\s{5}\s{4}\s{3}1$&$166934$ $901421$\\\hline
$2\s{3}\s{4}5\s{1}$&$166938$&$234\s{5}1$&$166939$&{\bf23451}&$166941$&$3\s{5}41\s{2}$&$166942$\\\hline
$35\s{4}12$&$166943$&$452\s{3}1$&$166945$&$2543\s{1}$&$166950$&$3254\s{1}$&$166951$\\\hline
$23\s{4}5\s{1}$&$166955$&$23\s{4}\s{5}1$&$166956$ $901718$&$2\s{3}4\s{5}1$&$166956$ $901724$&$23\s{4}51$&$166957$\\\hline
$2\s{3}\s{5}1\s{4}$&$166959$&$2\s{3}\s{5}\s{4}1$&$166969$&$2\s{5}\s{4}\s{1}\s{3}$&$166974$&$2\s{3}514$&$166978$\\\hline
$3\s{5}\s{4}\s{2}\s{1}$&$166980$&$243\s{5}\s{1}$&$166982$&$2\s{4}\s{3}51$&$166983$&$2354\s{1}$&$166985$ $921184$\\\hline
{\bf35421}&$166985$ $902215$&$234\s{5}\s{1}$&$166991$&$2\s{3}\s{4}51$&$166992$ $902202$&$3\s{5}\s{2}4\s{1}$&$166992$ $902120$\\\hline
$4\s{5}321$&$166992$ $902206$&$35\s{2}\s{4}1$&$166997$&$2435\s{1}$&$166998$ $902230$&$25\s{1}\s{4}3$&$166998$ $902155$\\\hline
$2\s{5}\s{1}4\s{3}$&$167001$&$254\s{1}3$&$167004$&$5\s{4}32\s{1}$&$167006$&$2345\s{1}$&$167008$\\\hline
$25\s{4}3\s{1}$&$167009$&$45\s{3}\s{2}1$&$167010$&$2\s{5}4\s{3}1$&$167011$&$45\s{2}\s{3}1$&$167014$\\\hline
$25\s{4}\s{1}3$&$167031$&$2\s{5}4\s{1}\s{3}$&$167034$&$2\s{4}153$&$167068$&$241\s{5}\s{3}$&$167091$\\\hline
$4\s{5}23\s{1}$&$167106$&$2\s{5}143$&$167110$&$251\s{4}\s{3}$&$167111$&$53\s{4}2\s{1}$&$167122$\\\hline
$35\s{4}1\s{2}$&$167131$&{\bf24153}&$167133$&$4\s{5}231$&$167139$&$3\s{4}512$&$167141$\\\hline
$235\s{1}4$&$167143$ $903551$&$34\s{5}12$&$167143$ $903656$&$23\s{5}\s{1}\s{4}$&$167144$&$4\s{5}\s{2}31$&$167158$\\\hline
$3\s{4}51\s{2}$&$167161$&$34\s{5}1\s{2}$&$167163$&$5\s{3}\s{4}2\s{1}$&$167188$&{\bf34512}&$167202$\\\hline
$4523\s{1}$&$167277$&$5342\s{1}$&$167300$&$3\s{5}412$&$167321$&$35\s{4}\s{2}1$&$167330$\\\hline
$3\s{5}42\s{1}$&$167332$&$13\s{2}5\s{4}$&$167408$&$153\s{4}\s{2}$&$167560$&&$167561$ $905557$\\[-1ex]
&&&&&&\rb{$2\s{1}4\s{5}\s{3}$}&$5067054$\\\hline
&$167561$ $905557$&$145\s{2}3$&$167601$&&$167602$ $906143$&&$167602$ $906143$\\[-1ex]
\rb{$2\s{1}\s{4}53$}&$5067055$&&&\rb{$2\s{1}453$}&$5073953$ $29335370$&\rb{$2\s{1}\s{4}\s{5}\s{3}$}&$5073953$ $29335426$\\\hline
$1452\s{3}$&$167646$&$35\s{1}42$&$167670$&$54\s{3}2\s{1}$&$167744$&&$167748$ $907383$\\[-1ex]
&&&&&&\rb{$2\s{1}4\s{5}3$}&$5083238$ $29397202$\\\hline
&$167748$ $907383$&$1534\s{2}$&$167749$ $907398$&$325\s{4}1$&$167749$ $907418$&$52\s{4}\s{3}\s{1}$&$167815$\\[-1ex]
\rb{$2\s{1}\s{4}5\s{3}$}&$5083238$ $29397203$&&&&&&\\\hline
&$167818$ $907708$&&$167818$ $907708$&$2\s{4}\s{5}\s{3}\s{1}$&$167826$&{\bf24531}&$167828$\\[-1ex]
\rb{$2\s{1}35\s{4}$}&$5083642$ $29380782$&\rb{$2\s{1}3\s{5}4$}&$5083642$ $29380784$&&&&\\\hline
$4\s{5}3\s{2}1$&$167832$&$5243\s{1}$&$167833$&$45\s{3}21$&$167835$&$13\s{5}\s{2}4$&$167844$\\\hline
$2\s{4}5\s{3}\s{1}$&$167848$&$24\s{5}31$&$167850$&$135\s{4}\s{2}$&$167855$ $908182$&$14\s{3}5\s{2}$&$167855$ $908181$\\\hline
\end{tabular}}\\
\vspace*{1ex}
\hfill{\tiny continued}
\end{table}

\begin{table}[ht]
{\tiny
\begin{tabular}{|cl|cl|cl|cl|}\hline
$14\s{3}\s{5}\s{2}$&$167863$&$351\s{4}\s{2}$&$167869$&$13\s{5}42$&$167877$&$3514\s{2}$&$167886$\\\hline
$32\s{5}4\s{1}$&$167923$&$32\s{5}41$&$167940$&$2\s{3}\s{5}4\s{1}$&$167942$ $909327$&$2\s{5}\s{4}3\s{1}$&$167942$ $909336$\\\hline
$235\s{4}1$&$167943$&$254\s{3}1$&$167944$&$2\s{4}15\s{3}$&$167951$&$23\s{5}4\s{1}$&$167959$\\\hline
$1453\s{2}$&$167960$ $909582$&$2\s{3}5\s{4}1$&$167960$ $909568$&$2\s{3}5\s{1}\s{4}$&$167961$&$2\s{3}\s{5}\s{1}4$&$167962$\\\hline
$2\s{4}53\s{1}$&$167963$&$24\s{5}\s{3}1$&$167965$&$23\s{5}14$&$167967$&$2351\s{4}$&$167968$ $909719$\\\hline
$25\s{3}14$&$167968$ $909740$&$2\s{4}513$&$167974$&$2\s{4}\s{5}1\s{3}$&$167977$&$523\s{4}\s{1}$&$167981$ $909851$\\\hline
$52\s{3}4\s{1}$&$167981$ $909855$&$2\s{5}\s{3}\s{1}4$&$167988$&$35\s{1}\s{4}2$&$167990$&$241\s{5}3$&$167991$\\\hline
$2\s{5}1\s{4}3$&$167993$&$14\s{5}\s{2}3$&$167998$ $910090$&$2531\s{4}$&$167998$ $910112$&$5234\s{1}$&$167998$ $910078$\\\hline
$45\s{3}1\s{2}$&$168007$&{\bf25314}&$168008$ $910322$&$3\s{5}\s{4}\s{2}1$&$168008$ $910269$&$4531\s{2}$&$168008$ $910276$\\\hline
$2514\s{3}$&$168011$ $910256$&$453\s{2}1$&$168011$ $910347$&$135\s{2}\s{4}$&$168012$&$4\s{5}3\s{2}\s{1}$&$168024$\\\hline
$24\s{5}3\s{1}$&$168027$&$2\s{4}5\s{3}1$&$168029$ $910494$&$3\s{5}2\s{4}1$&$168029$ $910481$&&$168039$
$909957$\\[-1ex]
&&&&&&\rb{$2\s{1}54\s{3}$}&$5104177$ $29555753$\\\hline
&$168039$ $909957$&$2\s{4}35\s{1}$&$168054$&$24\s{3}\s{5}1$&$168055$&$24\s{3}51$&$168056$\\[-1ex]
\rb{$2\s{1}\s{5}43$}&$5104177$ $29555755$&&&&&&\\\hline
$45\s{3}2\s{1}$&$168084$&&$168088$ $910579$&&$168088$ $910579$&$25\s{3}4\s{1}$&$168108$\\[-1ex]
&&\rb{$2\s{1}45\s{3}$}&$5110667$ $29617694$&\rb{$2\s{1}\s{4}\s{5}3$}&$5110667$ $29617699$&&\\\hline
$2\s{5}3\s{4}1$&$168109$&$25\s{4}\s{3}\s{1}$&$168116$&$2\s{5}431$&$168118$&$24\s{1}\s{5}\s{3}$&$168123$\\\hline
$325\s{4}\s{1}$&$168133$&$2\s{3}5\s{4}\s{1}$&$168134$&$23\s{5}41$&$168135$&$25\s{3}1\s{4}$&$168136$\\\hline
$3\s{5}\s{4}\s{1}\s{2}$&$168137$&$253\s{4}\s{1}$&$168140$&$2\s{5}\s{3}41$&$168141$&$2415\s{3}$&$168146$\\\hline
$24\s{5}\s{1}\s{3}$&$168147$ $911472$&$2\s{5}\s{3}\s{4}1$&$168147$ $911476$&$354\s{1}2$&$168152$&$35\s{4}21$&$168155$\\\hline
$2\s{3}\s{5}14$&$168159$&$2\s{3}51\s{4}$&$168160$ $911630$&$25\s{4}1\s{3}$&$168160$ $911639$&$2534\s{1}$&$168163$ $911669$\\\hline
$4532\s{1}$&$168163$ $911687$&$2\s{4}\s{5}3\s{1}$&$168166$&$245\s{1}3$&$168167$&$245\s{3}1$&$168168$ $911687$\\\hline
$25\s{3}\s{1}4$&$168168$ $911692$&$235\s{4}\s{1}$&$168169$&$2\s{3}\s{5}41$&$168170$ $911718$&$3\s{5}\s{4}2\s{1}$&$168170$ $911823$\\\hline
$24\s{5}\s{3}\s{1}$&$168174$&$2\s{4}531$&$168176$&$2\s{5}31\s{4}$&$168177$&$3\s{5}4\s{2}\s{1}$&$168184$\\\hline
$2\s{4}\s{1}5\s{3}$&$168200$&$154\s{3}\s{2}$&$168202$&$354\s{2}1$&$168203$&$2541\s{3}$&$168207$\\\hline
$24\s{5}13$&$168211$&$2\s{4}\s{5}\s{3}1$&$168212$&$3\s{5}4\s{2}1$&$168215$&$35\s{4}2\s{1}$&$168216$\\\hline
$3542\s{1}$&$168217$&$3\s{5}241$&$168219$&$2453\s{1}$&$168228$&$3524\s{1}$&$168255$\\\hline
$2\s{4}51\s{3}$&$168265$&$145\s{3}\s{2}$&$168266$&$3\s{2}5\s{4}1$&$168268$&$24\s{3}5\s{1}$&$168279$\\\hline
$24\s{3}\s{5}\s{1}$&$168280$&$2\s{4}351$&$168281$&$2\s{5}\s{1}43$&$168292$&$25\s{3}\s{4}\s{1}$&$168296$\\\hline
$2\s{5}341$&$168297$&$3\s{4}5\s{2}1$&$168300$&$2\s{5}3\s{1}4$&$168304$ $912844$&{\bf34521}&$168304$ $913052$\\\hline
$25\s{1}\s{4}\s{3}$&$168308$ $912905$&$52\s{4}3\s{1}$&$168308$ $912922$&$35\s{4}\s{1}2$&$168312$&$34\s{5}21$&$168317$ $913171$\\\hline
$34\s{5}2\s{1}$&$168317$ $913172$&$23\s{5}\s{1}4$&$168328$ $913181$&$3\s{5}4\s{1}\s{2}$&$168328$ $913277$&$145\s{2}\s{3}$&$168330$ $913130$\\\hline
$3\s{4}5\s{2}\s{1}$&$168330$ $913304$&$2\s{5}4\s{1}3$&$168333$&$3541\s{2}$&$168343$&$235\s{1}\s{4}$&$168344$\\\hline
$3\s{4}52\s{1}$&$168353$&$253\s{1}\s{4}$&$168354$&$24\s{1}\s{5}3$&$168355$&$3\s{5}\s{2}\s{4}1$&$168361$\\\hline
$253\s{1}4$&$168363$ $913662$&$2\s{5}\s{4}\s{1}3$&$168363$ $913651$&$2451\s{3}$&$168366$&$2\s{4}5\s{1}\s{3}$&$168367$\\\hline
$34\s{5}\s{2}1$&$168369$&$2\s{5}413$&$168386$&$3\s{4}521$&$168389$&$352\s{4}\s{1}$&$168394$\\\hline
$4\s{5}3\s{1}2$&$168396$&$2\s{5}\s{4}13$&$168397$&$345\s{2}1$&$168402$&$3\s{5}\s{4}21$&$168423$\\\hline
$3\s{5}\s{4}\s{1}2$&$168431$&$2\s{4}\s{5}13$&$168435$ $914602$&$34\s{5}\s{2}\s{1}$&$168435$ $914677$&$24\s{5}\s{1}3$&$168438$\\\hline
$3\s{2}\s{5}41$&$168460$&$5\s{3}\s{4}\s{2}\s{1}$&$168475$&$53\s{4}\s{2}\s{1}$&$168486$&$3451\s{2}$&$168493$\\\hline
$3\s{4}5\s{1}\s{2}$&$168509$&$3\s{5}\s{4}12$&$168515$&$3\s{5}\s{2}41$&$168521$&$2\s{4}\s{5}\s{1}3$&$168522$\\\hline
$3452\s{1}$&$168525$&$25\s{1}4\s{3}$&$168526$&$24\s{1}5\s{3}$&$168527$ $915136$&$2\s{5}\s{1}\s{4}3$&$168527$ $915161$\\\hline
$34\s{5}\s{1}\s{2}$&$168527$ $915307$&$3\s{5}4\s{1}2$&$168537$&$25\s{4}\s{1}\s{3}$&$168542$&$254\s{3}\s{1}$&$168546$\\\hline
$2\s{5}\s{4}31$&$168547$&$3\s{5}421$&$168554$&$34\s{5}\s{1}2$&$168563$&$35\s{2}4\s{1}$&$168567$\\\hline
$35\s{4}\s{2}\s{1}$&$168583$&$245\s{3}\s{1}$&$168584$&$2\s{4}\s{5}31$&$168585$&$254\s{1}\s{3}$&$168587$\\\hline
$543\s{2}\s{1}$&$168588$&$45\s{3}\s{2}\s{1}$&$168597$&$3\s{5}\s{1}42$&$168621$&$245\s{1}\s{3}$&$168625$\\\hline
$35\s{1}4\s{2}$&$168636$&$452\s{3}\s{1}$&$168648$&$354\s{2}\s{1}$&$168661$&$3\s{2}5\s{4}\s{1}$&$168670$\\\hline
$354\s{1}\s{2}$&$168670$&$345\s{2}\s{1}$&$168673$&$345\s{1}\s{2}$&$168682$&$524\s{3}\s{1}$&$168691$\\\hline
$35\s{4}\s{1}\s{2}$&$168745$&$534\s{2}\s{1}$&$168757$&$35\s{2}\s{4}\s{1}$&$168760$&$45\s{2}\s{3}\s{1}$&$168766$\\\hline
$453\s{2}\s{1}$&$168820$&$453\s{1}\s{2}$&$168829$&&&&\\\hline
\end{tabular}
\vspace*{1ex}
\caption{
I-Wilf classes of $B_5$ and the numbers $|SI_n(\tau)|$ for $n=9,10,11,12$.
To determine the class to which the pattern $\s{1}\s{4}\s{5}23$ belongs, calculate $|SI_9(\s{1}\s{4}\s{5}23)|=168330$.
This number corresponds to both the patterns $145\s{2}\s{3}$ and $3\s{4}5\s{2}\s{1}$ above. 
To decide which of these is the correct one, it is necessary to calculate $|SI_{10}(\s{1}\s{4}\s{5}23)|=913130$. Thus $\s{1}\s{4}\s{5}23$ belongs to the class represented by $145\s{2}\s{3}$.
}
\vspace*{-2ex}}
\end{table}

\begin{table}[ht]
{\tiny
\begin{tabular}{|cl|cl|cl|cl|cl|cl|}\hline
$361542$&$97405 $&$465132$&$97511 $&$361452$&$98805 $&$351624$&$99133 $&$426153$&$99287 $&$146253$&$99321 $\\ \hline
$132546$&$99432 $&$125436$&$99521 $&$154326$&$99585 $&$153624$&$99650 $&$124356$&$99653 $&$123546$&$99729 $\\ \hline
$624351$&$99857 $&$625431$&$99885 $&$123456$&$99991 $&$623541$&$100021$&$645231$&$100088$&$632541$&$100156$\\ \hline
$563412$&$100293$&$623451$&$100615$&$163542$&$100879$&$463152$&$100992$&$164352$&$101197$&$125634$&$101405$\\ \hline
$156423$&$101451$&$145236$&$101662$&$126453$&$101754$&$163452$&$101918$&$153426$&$102109$&$135426$&$104236$\\ \hline
$136542$&$105312$&$124653$&$105971$&$124536$&$106788$&$154362$&$106857$&$156342$&$107185$&$125463$&$107578$\\ \hline
$326154$&$107772$&$134526$&$108083$&$136254$&$108336$&$265431$&$108967$&$143625$&$108969$&$145326$&$109293$\\ \hline
$261543$&$109404$&$143652$&$109443$&$462513$&$109514$&$132564$&$109674$&$135246$&$109943$&$136452$&$110137$\\ \hline
$123564$&$110264$&$134652$&$110707$&$124563$&$110872$&$135462$&$110964$&$146352$&$111024$&$143562$&$111229$\\ \hline
$635421$&$111594$&$264351$&$111647$&$135624$&$111648$&$263541$&$111733$&$153462$&$111836$&$124635$&$111871$\\ \hline
$362541$&$111963$&$125643$&$112058$&$624531$&$112186$&$462531$&$112231$&$156432$&$112493$&$261453$&$112598$\\ \hline
$153642$&$112738$&$253614$&$112805$&$145263$&$112830$&$246153$&$112962$&$134625$&$113031$&$326541$&$113101$\\ \hline
$134562$&$113121$&$463251$&$113154$&$236154$&$113168$&$263451$&$113331$&$362451$&$113424$&$164532$&$113439$\\ \hline
$154623$&$113690$&$136524$&$113837$&$426513$&$113909$&$136245$&$114046$&$351642$&$114060$&$236541$&$114071$\\ \hline
$254361$&$114129$&$462351$&$114245$&$146325$&$114470$&$256341$&$114598$&$326514$&$114730$&$146523$&$114833$\\ \hline
$146532$&$115050$&$364152$&$115051$&$562431$&$115131$&$251634$&$115165$&$463512$&$115289$&$564321$&$115297$\\ \hline
$261354$&$115305$&$243615$&$115357$&$264513$&$115506$&$365142$&$115532$&$324651$&$115600$&$635241$&$115605$\\ \hline
$256413$&$115714$&$243651$&$115741$&$264153$&$115762$&$634521$&$116018$&$564231$&$116084$&$154632$&$116098$\\ \hline
$264531$&$116206$&$365421$&$116214$&$265413$&$116546$&$241653$&$116580$&$234651$&$116603$&$135642$&$116656$\\ \hline
$145362$&$116665$&$562341$&$116676$&$236514$&$116688$&$235461$&$116747$&$251364$&$117002$&$645321$&$117190$\\ \hline
$465312$&$117342$&$234615$&$117530$&$135264$&$117649$&$234561$&$117661$&$325614$&$117792$&$256314$&$118369$\\ \hline
$265143$&$118372$&$231564$&$118450$&$231645$&$118517$&$346152$&$118533$&$563421$&$118646$&$326451$&$118724$\\ \hline
$145623$&$118881$&$465321$&$119049$&$264315$&$119084$&$246513$&$119204$&$136425$&$119269$&$251643$&$119284$\\ \hline
$236145$&$119306$&$261534$&$119411$&$256431$&$119481$&$426531$&$119592$&$256134$&$119745$&$236451$&$119864$\\ \hline
$456312$&$120024$&$356412$&$120049$&$356142$&$120195$&$364251$&$120269$&$235614$&$120277$&$254613$&$120434$\\ \hline
$265341$&$120451$&$362514$&$120655$&$253461$&$120790$&$246351$&$120922$&$254631$&$121026$&$365412$&$121073$\\ \hline
$246315$&$121125$&$465231$&$121289$&$263154$&$121348$&$145632$&$121395$&$263514$&$121571$&$251463$&$121692$\\ \hline
$254163$&$121697$&$235164$&$121719$&$253641$&$121786$&$263415$&$121892$&$325641$&$121936$&$246135$&$121959$\\ \hline
$246531$&$122125$&$356241$&$122422$&$245163$&$122425$&$426351$&$122452$&$256143$&$122484$&$436512$&$122608$\\ \hline
$241635$&$122668$&$364521$&$122725$&$352641$&$122840$&$235641$&$122894$&$245613$&$122957$&$245361$&$123195$\\ \hline
$346251$&$123251$&$463521$&$123375$&$465213$&$123413$&$456132$&$123474$&$364512$&$123518$&$456231$&$123756$\\ \hline
$236415$&$123833$&$356214$&$123835$&$354621$&$123935$&$365241$&$124192$&$346512$&$124405$&$356124$&$124936$\\ \hline
$265134$&$125054$&$265314$&$125541$&$245631$&$125665$&$365214$&$125736$&$356421$&$126250$&$345612$&$126268$\\ \hline
$436521$&$126552$&$346521$&$126743$&$354612$&$127013$&$456321$&$127598$&$345621$&$128803$&
& \\\hline
\end{tabular}
\vspace*{1ex}
\caption{I-Wilf classes of $S_6$ and the numbers $|I_{12}(\tau)|$}
\vspace*{-2ex}}
\end{table}

\end{document}